\documentclass{amsart}
\usepackage{amsmath,amssymb}
\newtheorem{theorem}{Theorem}
\newtheorem{lemma}[theorem]{Lemma}
\newtheorem{proposition}[theorem]{Proposition}

\newtheorem{remark}[theorem]{Remark}
\newtheorem{example}[theorem]{Example}

\begin{document}
\title[separable products of free modules]{Valuation domains whose products of free modules are separable}
\author{Fran\c cois Couchot}
\address{Laboratoire de Math\'ematiques Nicolas Oresme, CNRS UMR
  6139,
D\'epartement de math\'ematiques et m\'ecanique,
14032 Caen cedex, France}
\email{couchot@math.unicaen.fr} 



\begin{abstract} It is proved that if $R$ is a valuation domain with maximal ideal $P$ and if $R_L$ is countably generated for each prime ideal $L$,  then $R^R$ is separable if and only $R_J$ is maximal, where $J=\cap_{n\in\mathbb{N}}P^n$.
\end{abstract}
\maketitle

When $R$ is a valuation domain satisfying one of the following two conditions:
\begin{enumerate}
\item $R$ is almost maximal and its quotient field $Q$ is countably generated
\item $R$ is archimedean
\end{enumerate}
Franzen proved in \cite{Fra84} that $R^{\mathbb{N}}$ is separable if and only if $R$ is maximal or discrete of rank one. In \cite[Theorem XVI.5.4]{FuSa01}, Fuchs and Salce gave a slight generalization of this result and showed that $R^{\mathbb{N}}$ is separable if and only if $R$ is  discrete of rank one, when $R$ is slender. The aim of this paper is to give another generalization of Franzen's result by proving Theorem~\ref{T:main} below. If the maximal ideal $P$ is principal, we get that $R^R$ can be separable when $R$ is neither maximal nor discrete of rank one. This is a negative answer to \cite[Problem 59]{FuSa01}. For proving his result, Franzen began by showing that each  archimedean valuation domain which is not almost maximal, possesses an indecomposable reflexive module of rank 2. We  use a similar argument in the proof of Theorem~\ref{T:main}. Finally we  give an example of a non-archimedean non-slender valuation domain such that $R^{\mathbb{N}}$ is not separable. This is a positive answer to \cite[Problem 58]{FuSa01}.

\medskip

In the sequel, $R$ is a commutative unitary ring.  An $R$-module whose submodules are totally ordered by inclusion, is said to be \textbf{uniserial}. If $R$ is a uniserial $R$-module, we say that $R$ is a \textbf{valuation ring}. 

The \textbf{$R$-topology} of $R$ is the linear topology for which each non-zero
ideal is a neighborhood of $0$. When $R$ is a valuation ring with maximal ideal $P$ and $A$ is a proper ideal, then $R/A$ is Hausdorff in the $R/A$-topology if and only if $A\ne Pa,\ \forall 0\ne a\in R$. We say that $R$ is \textbf{(almost) maximal}  if $R/A$ is complete in the $R/A$-topology for each (non-zero) proper ideal $A\ne Pa,\ \forall 0\ne a\in R$.

From now on, $R$ is a valuation domain, $P$ is its maximal ideal and $Q$ is its field of quotients. Let $M$ be an $R$-module and let $N$ be a submodule. We say that $N$ is a \textbf{pure submodule} of $M$ if $rN=rM\cap N,\ \forall r\in R$.
Let $M$ be a torsion-free module. We say that $M$ is \textbf{separable} if each pure uniserial submodule is a summand. Recall that each element $x$ of $M$ is contained in a pure uniserial submodule $U$, where  $U$ is the inverse image of the torsion submodule of $M/Rx$ by the canonical map $M\rightarrow M/Rx$. 
Let $M$ be a non-zero $R$-module. As in
\cite{FuSa01} we set:
 \[M^{\sharp}=\{s\in R\mid sM\subset M\}.\]
Then $M^{\sharp}$ is a prime ideal. We say that an ideal $A$ is \textbf{archimedean} if $A^{\sharp}=P$.

\begin{proposition} \label{P:sepa}  Let $J=\cap_{n\in\mathbb{N}}P^n$. Then,  $R^{\Lambda}$ is separable for each index set $\Lambda$ if $R_J$ is maximal.
\end{proposition}

\textbf{Proof.}  If $P$ is not finitely generated then $J=P$. In this case $R$ is maximal, whence $R^{\Lambda}$ is separable by \cite[Theorem 51]{Mat72}. Suppose now that $P=Rp$ for some $p\in P$. Let $U$ be a pure uniserial submodule of $R^{\Lambda}$. We must prove that $U$ is a summand. First assume that  $U^{\sharp}=P$. Then $pU\ne U$, whence $U=Ru$ for some $u\in U\setminus pU$. If $u=(u_{\lambda})_{{\lambda}\in{\Lambda}}$, there exists $\mu\in{\Lambda}$ such that $u_{\mu}$ is a unit because $pU=U\cap pR^{\Lambda}$. Then in the product $R^{\Lambda}$, the $\mu$th component can be replaced by $Ru$. So, $U$ is a summand. Now assume that $U^{\sharp}\subseteq J$. It follows that $U$ is a pure uniserial $R_J$-submodule of $(R^{\Lambda})_J$. Since $R_J$ is maximal, $U$ is a summand of $(R^{\Lambda})_J$. Then $U$ is a summand of $R^{\Lambda}$ too.\qed 

\bigskip
From Proposition~\ref{P:sepa} we deduce the following example which gives a negative answer to \cite[Problem 59]{FuSa01}.
\begin{example} \label{E:59} \textnormal{Let $R=\mathbb{Z}_p+X\mathbb{Q}[[X]]$, where $p$ is a prime number and $\mathbb{Z}_p$ is the localization of $\mathbb{Z}$ at the prime ideal $p\mathbb{Z}$. Then $J=X\mathbb{Q}[[X]]$, $R/J\cong\mathbb{Z}_p$ and $R_J\cong\mathbb{Q}[[X]]$. It follows that $R$ is neither maximal nor discrete of rank one, but $R_J$ is maximal, whence $R^{\Lambda}$ is separable for each index set $\Lambda$ by Proposition~\ref{P:sepa}. So, \cite[Exercise XVI.5.5]{FuSa01} is wrong.}
\end{example}

\bigskip
To prove Theorem~\ref{T:main} some preliminary results are needed.

If $M$ is an $R$-module, $\mathrm{Hom}_R(M,R)$ is denoted by $M^*$ and $\lambda_M:M\rightarrow M^{**}$ is the canonical map. We say that $M$ is \textbf{reflexive} if $\lambda_M$ is an isomorphism.
An $R$-module $F$ is \textbf{pure-injective} if for every pure exact
sequence
$0\rightarrow N\rightarrow M\rightarrow L\rightarrow 0$
 of $R$-modules, the following sequence
 \[0\rightarrow\mathrm{Hom}_R(L,F)
\rightarrow\mathrm{Hom}_R(M,F)\rightarrow\mathrm{Hom}_R(N,F)\rightarrow
0\] is exact. An $R$-module $B$
is a \textbf{pure-essential extension} of a submodule $A$ if $A$ is a
pure submodule of $B$ and, if for each
submodule $K$ of $B$, either $K\cap A\ne 0$ or $(A+K)/K$ is not a pure
submodule of $B/K$. 
We say that $B$ is a \textbf{pure-injective hull}
of $A$ if $B$ is pure-injective and a pure-essential extension of $A$.
 By \cite[Chapter XIII]{FuSa01} each $R$-module $M$ has a
 pure-injective hull and any two pure-injective hulls of $M$ are isomorphic.
For any module $M$, we denote by $\widehat{M}$ its pure-injective hull. If $S$ is a maximal immediate extension of $R$, then $S\cong\widehat{R}$ by \cite[Proposition XIII.5.1]{FuSa01}. For each $s\in S\setminus R$, $\mathrm{B}(s)=\{r\in R\mid s\notin R+rS\}$ is called the \textbf{breadth ideal} of $s$.

\begin{proposition} \label{P:indecom}
Let $A$ be a non-zero archimedean ideal such that $A\ne Pa$ for each $a\in R$. Assume that $R/A$ is not complete in the $R/A$-topology. Then there exists an indecomposable reflexive module of rank $2$.
\end{proposition}
\textbf{Proof.}  Since $R/A$ is not complete in the $R/A$-topology, by \cite[Lemma V.6.1]{FuSa01} there exists $x\in\widehat{R}\setminus R$ such that $A=\mathrm{B}(x)$.  Let $U$ be a pure uniserial submodule of $\widehat{R}/R$ containing $x+R$ and let $M$ be the inverse image of $U$ by the natural map $\widehat{R}\rightarrow\widehat{R}/R$. Then  $M$ is a pure submodule of $\widehat{R}$. By \cite[Example XV.6.1]{FuSa01} $M$ is indecomposable.  Since $M$ is a pure extension of $R$ by $U$ then $\mathrm{Ext}^1_R(U,R)\ne 0$ and $U$ is torsion-free.  Now, we show that $U^{\sharp}=P$. Let $0\ne s\in P$. Then $A\subset s^{-1}A$. Let $t\in (R\cap s^{-1}A)\setminus A$. Therefore $x=r+ty$ for some $r\in R$ and $y\in\widehat{R}$. Since $M$ is a pure submodule of $\widehat{R}$, we may assume that $y\in M$. By \cite[Lemma 1.3]{SaZa85}, $\mathrm{B}(y)=t^{-1}A$. Consequently $y+R\notin sU$ and $U^{\sharp}=P$. Since $U$ is a torsion-free module of rank one and $U^{\sharp}\ne 0$, $U$ is isomorphic to a proper submodule of $Q$. So $U$ is isomorphic to an ideal of $R$. By \cite[Proposition 3.3]{Fra84}, the proof is complete.
\qed

\begin{proposition} \label{P:locsep} Assume that $R^{\Lambda}$ is separable for an index set $\Lambda$. Then $(R_L)^{\Lambda}$ is separable for each prime ideal $L$.
\end{proposition}
\textbf{Proof.} The assertion is obvious if $L=0$. Now suppose $L\ne 0$ and let $0\ne a\in L$. Then $R_La$ is an ideal contained in $L$. Let $U$ be a pure uniserial submodule of $(aR_L)^{\Lambda}$ and let $V$ be the inverse image of the torsion submodule of $R^{\Lambda}/U$ by the surjection of $R^{\Lambda}$ onto $R^{\Lambda}/U$. Then $V$ is a pure uniserial submodule of $R^{\Lambda}$. Let $p$ be a projection of $R^{\Lambda}$ onto $V$ and $q=p\vert_{(aR_L)^{\Lambda}}$. For each $s\in R\setminus aR_L$ we have $(aR_L)^{\Lambda}\subseteq sR^{\Lambda}$. Thus $\mathrm{Im}\ q\subseteq sR^{\Lambda}$. Since $aR_L=\cap_{s\in R\setminus aR_L}sR$ we get $\mathrm{Im}\ q\subseteq (aR_L)^{\Lambda}$. On the other hand $U\subseteq\mathrm{Im}\ q$ and the equality holds because $U$ is a pure submodule. It is obvious that $(aR_L)^{\Lambda}=a(R_L)^{\Lambda}\cong (R_L)^{\Lambda}$. Hence $(R_L)^{\Lambda}$ is separable.\qed

\begin{lemma} 
\label{L:locaComp} Let $L$ be a prime ideal of $R$ and let $A$ be a proper ideal of $R_L$. If $R/A$ is complete in the $R/A$-topology then $R_L/A$ is also complete in the $R_L/A$-topology.
\end{lemma}
\textbf{Proof.} Let $(a_i+A_i)_{i\in I}$ be a family of cosets of $R_L$ such that $a_i\in a_j+A_j$ if $A_i\subset A_j$ and such that $A=\cap_{i\in I}A_i$. We may assume that $A_i\subseteq L,\ \forall i\in I$. So, $a_i+L=a_j+L,\ \forall i,j\in I$. Let $b\in a_i+L,\ \forall i\in I$. It follows that $a_i-b\in L,\ \forall i\in I$. Since $R/A$ is complete in the $R/A$-topology, $\exists c\in R$ such that $c+b-a_i\in A_i,\ \forall i\in I$. Hence $R_L/A$ is complete in the $R_L/A$-topology too. \qed

\medskip
Recall that an $R$-module $M$ is \textbf{slender} if for every morphism $f:R^{\mathbb{N}}\rightarrow M$, there exists $n_0\in\mathbb{N}$ such that $f(e_n)=0,\ \forall n\geq n_0$, where $e_n=(\delta_{n,m})_{m\in\mathbb{N}}$. In the proof of Theorem~\ref{T:main} we will use the following result: 
\begin{proposition} \label{P:slender} \cite[Corollary 21]{Dim83} Let $R$ be a valuation domain such that $Q$ is countably generated. Then $R$ is slender if and only if $R$ is not complete in the $R$-topology.
\end{proposition}

\medskip
 The following proposition can be easily proved.
\begin{proposition} \label{P:countable}
The following conditions are equivalent:
\begin{enumerate}
\item $R_L$ is countably generated for each prime ideal $L$.
\item For each prime ideal $L$ which
  is the intersection of the
  set of primes containing properly $L$ there is a countable subset
  whose intersection is $L.$
\item For each prime ideal $L$, the quotient field of $R/L$ is countably generated.
\end{enumerate}
\end{proposition}
\begin{theorem} \label{T:main}
Assume that $R$ satisfies the equivalent conditions of Proposition~\ref{P:countable}. Let $J=\cap_{n\in\mathbb{N}}P^n$. Then the following conditions are equivalent:
\begin{enumerate}
\item $R^{\Lambda}$ is separable for each index set $\Lambda$;
\item $R^R$ is separable;
\item $R_J$ is maximal.
\end{enumerate}
Moreover, if each ideal is countably generated then these conditions are equivalent to: 
$R^{\mathbb{N}}$ is separable.
\end{theorem}
\textbf{Proof.} It is obvious that $(1)\Rightarrow (2)$. By Proposition~\ref{P:sepa}, $(3)\Rightarrow (1)$.

$(2)\Rightarrow (3)$. We must prove that $R_J/A$ is complete in the $R_J/A$-topology for each ideal $A$ of $R_J$, where $A\ne Jr$, $\forall 0\ne r\in R$. By Lemma~\ref{L:locaComp} it is enough to show that $R/A$ is complete in the $R/A$-topology.

First we assume that $A$ is prime, $A\subset J$. Suppose that $R/A$ is not complete in the $R/A$-topology. By \cite[Lemma XVI.5.3]{FuSa01}, $(R/A)^{\mathbb{N}}$ is separable. Since $R$ satisfies the conditions of Proposition~\ref{P:countable}, the quotient field of $R/A$ is countably generated. It follows by Proposition~\ref{P:slender} that $R/A$ is slender.  By \cite[Theorem XVI.5.4]{FuSa01} $R/A$ is a discrete valuation domain of rank one . Clearly we get a contradiction. Hence $R/A$ is complete in the $R/A$-topology. Suppose that $A=rA^{\sharp}$ for some $0\ne r\in R$ where $A^{\sharp}\subset J$. It is easy to deduce  the completeness of $R/A$ from that of $R/A^{\sharp}$.

Now assume that $A\ne rA^{\sharp}$, $\forall r\in R$. First we show that $R_{A^{\sharp}}/A$ is complete in the $R_{A^{\sharp}}/A$-topology.  By way of contradiction, suppose it is not true. We put $R'=R_{A^{\sharp}}$ and $N^*=\mathrm{Hom}_{R'}(N,R')$ if $N$ is an $R'$-module. Then $A$ is an archimedean ideal of $R'$. By Proposition~\ref{P:indecom} there exists an indecomposable reflexive $R'$-module $M$ of rank $2$. The map $\phi:M^{**}\rightarrow (R')^{M^*}$ defined by $\phi(u)=(u(m))_{m\in M^*},\ \forall u\in M^{**}$, is a pure monomorphism. Since $M^*$ has the same cardinal as $R$, $(R')^{M^*}$ is separable by Proposition~\ref{P:locsep}. It follows that $M$ is separable. This contradicts that $M$ is indecomposable.

 Now we prove that $R/A$ is complete in the $R/A$-topology. Let $(a_i+A_i)_{i\in I}$ be a family of cosets of $R$ such that $a_i\in a_j+A_j$ if $A_i\subset A_j$ and such that $A=\cap_{i\in I}A_i$. We may assume that $A\subset A_i\subseteq A^{\sharp},\ \forall i\in I$. We put $A_i'=(A_i)_{A^{\sharp}},\ \forall i\in I$. We know that $A=\cap_{a\notin A}{A^{\sharp}}a$. Consequently, if $a\notin A$, there exists $i\in I$ such that  $A_i\subseteq A^{\sharp}a$, whence $A_i'\subseteq A^{\sharp}a$. It follows that $A=\cap_{i\in I}A_i'$. Clearly, $a_i\in a_j+A_j'$ if $A_i'\subset A_j'$. Then there exists $c\in R_{A^{\sharp}}$ such that $c\in a_i+A_i',\ \forall i\in I$. Since $A_i'\subset R,\ \forall i\in I$, $c\in R$. From $A=\cap_{j\in I}A_j'$ and $A\subset A_i,\ \forall i\in I$ we deduce that $\forall i\in I,\ \exists j\in I$ such that $A_j'\subset A_i$. We get that $c\in a_i+A_i$ because $c-a_j\in A_j'\subseteq A_i$ and $a_j-a_i\in A_i$. So, $R/A$ is complete in the $R/A$-topology.

To prove the last assertion it is enough to observe that $M^*$ is countably generated over $R'$ and consequently $M^{**}$ is isomorphic to a pure $R'$-submodule of $(R')^{\mathbb{N}}$. \qed

\begin{remark} \textnormal{ In proving that $R/A$ is complete, we use the hypothesis that $R$ satisfies the conditions of Proposition~\ref{P:countable} only when $A$ is isomorphic to a prime ideal. In the other case, this result can be obtained with the sole hypothesis that $R^R$ is separable.}
\end{remark}

\medskip
 So, even if $R$ doesn't satisfy the conditions of Proposition~\ref{P:countable} the next proposition holds:
\begin{proposition}
Let the notations be as in Theorem~\ref{T:main} and suppose that $R^R$ is separable. Then $R$ satisfies the following conditions:
\begin{enumerate}
\item $R/L$ is not slender for each prime ideal $L\subset J$.
\item $R/A$ is complete in the $R/A$-topology for each ideal $A$ which is not isomorphic to a prime ideal and such that $A^{\sharp}\subseteq J$.
\end{enumerate} 
\end{proposition}


The following example gives a positive answer to \cite[Problem 58]{FuSa01}.
\begin{example} 
\textnormal{Let $T$ be a non-discrete archimedean valuation domain which is not complete in the $T$-topology, let $K$ be its quotient field and let $R=T+XK[[X]]$. Let $L=XK[[X]]$. Then $Q$ and $R_L$ are countably generated. Moreover $R$ is complete in the $R$-topology because $R_L(\cong K[[X]])$ is maximal and  $R/L(\cong T)$ is not complete in the $R/L$-topology. So, $R$ is non-archimedean, $R^{\mathbb{N}}$ is not separable and $R$ is not slender.}
\end{example}


\begin{thebibliography}{1}

\bibitem{Dim83}
R.~Dimitri\'c.
\newblock Slender modules over domains.
\newblock {\em Comm. Algebra}, 11(15):1685--1700, (1983).

\bibitem{Fra84}
B.~Franzen.
\newblock On the separability of a direct product of free modules over a
  valuation domain.
\newblock {\em Arch. Math.}, 42:131--135, (1984).

\bibitem{FuSa01}
L.~Fuchs and L.~Salce.
\newblock {\em Modules over Non-Noetherian Domains}.
\newblock Number~84 in Mathematical Surveys and Monographs. American
  Mathematical Society, Providence, (2001).

\bibitem{Mat72}
E.~Matlis.
\newblock {\em Torsion-free modules}.
\newblock University of Chicago Press, (1972).

\bibitem{SaZa85}
L.~Salce and P.~Zanardo.
\newblock Some cardinals invariants for valuation domains.
\newblock {\em Rend. Sem. Mat. Univ. Padova}, 74:205--217, (1985).

\end{thebibliography}
\end{document}